\documentclass[12pt]{article}
\def\q{\hfill\rule{1ex}{1ex}}
\def\0{\emptyset}

\def\q{\hfill\rule{1ex}{1ex}}
\def\QEDopen{{\hfill\setlength{\fboxsep}{0pt}\setlength{\fboxrule}{0.2pt}\fbox{\rule[0pt]{0pt}{1.3ex}\rule[0pt]{1.3ex}{0pt}}}}

\usepackage{graphicx}
\usepackage{amsmath}
\usepackage{url}
\usepackage{graphicx}   
\usepackage{dsfont}
\usepackage{amsfonts}
\usepackage{amsmath}
\usepackage{bm}      
\usepackage{multirow,bigstrut,booktabs}    
\usepackage{setspace}
\usepackage{indentfirst}
\usepackage{latexsym,  amsthm, graphicx, mathrsfs} 
\usepackage{amssymb}
\usepackage{amsmath,amssymb,latexsym,color}
\usepackage{authblk}
\usepackage{url} 
\usepackage[mathscr]{eucal}
\usepackage[colorlinks,
            linkcolor=blue,
            anchorcolor=green,
            citecolor=magenta
           ]{hyperref}
\usepackage{geometry}
\begin{document}
\title{\bf Edge pancyclic derangement graphs
\thanks{
This work is  supported by the National Natural Science Foundation of China (Grant 12171272 \& 11971158) and  Tsinghua University Initiative Scientific Research Program}}
\author{{\small\bf Zequn
		Lv}\thanks{Email:  lvzq19@mails.tsinghua.edu.cn
	}\quad {\small\bf Mengyu Cao}\thanks{Corresponding author. Email: mcao@tsinghua.edu.cn}\quad {\small\bf Mei
		Lu}\thanks{Email: lumei@tsinghua.edu.cn}\\
	{\small Department of Mathematical Sciences, Tsinghua
		University, Beijing 100084, China.}}

\date{}

\maketitle\baselineskip 16.3pt

\begin{abstract}

We consider the derangement graph in which the vertices are permutations of $\{ 1,\ldots, n\}$. Two
vertices are joined by an edge if the corresponding permutations differ in every position. The
derangement graph is known to be Hamiltonian and Hamilton-connected.
In this note, we show that the derangement graph is edge pancyclic if $n\ge 4$.

\end {abstract}

\vskip.3cm

\section{Introduction}

Let $\Gamma=(V,E)$ be a graph. For a subset $S$ of $V$, $\Gamma[S]$ is a subgraph of $\Gamma$ induced by $S$. For any $v\in V(\Gamma)$, let $N(v)= \{w~|~ vw\in E(\Gamma)\}$ be the neighborhood of $v$ and $d(v)=|N(v)|$ be the degree of $v$. Let  $\delta=\delta(\Gamma)$ denote the  minimum degree of $\Gamma$. A matching
of size $s$ in $\Gamma$ is a family of $s$ pairwise disjoint edges of $\Gamma$. If the matching covers all
the vertices of $\Gamma$, then we call it a {\em perfect matching}. Given a graph $\Gamma$ of order $n \ge3$, we say that $\Gamma$ is {\em Hamiltonian}
if $\Gamma$ contains a cycle that spans $V(\Gamma)$, the vertex set of $\Gamma$. We say that $\Gamma$ is {\em pancyclic}
if $\Gamma$ contains a cycle of each possible length from 3 to $n$.  The graph $\Gamma$ is called {\em vertex pancyclic} ({\em edge pancyclic}) if every vertex (edge) is contained on a cycle of each possible length from
3 to $n$. Obviously, if $\Gamma$ is edge pancyclic then it is vertex pancyclic; if $\Gamma$ is vertex  pancyclic then it is pancyclic; if $\Gamma$ is  pancyclic then it is Hamiltonian.

Let $G$ be a finite group and $S$ be an inverse closed subset of $G$ with $1\notin S$. The {\em Cayley graph} $\Gamma (G, S)$ is the graph which has the elements of
$G$ as its vertices and two vertices $u, v \in G$ are joined by an edge if and only if $v = su$
for some $s \in S$. The graph $\Gamma (G, S)$ is
connected if and only if $S$ is a set of generators for $G$ and it is vertex-transitive.

Let $S_n$ be the symmetric group on $[n] = \{1,\ldots, n\}$. We let $D_n$ be the set of all derangements of $S_n$, where a derangement is a permutation without fixed points. Then $$|D_n|=n!\sum_{i=0}^n\frac{(-1)^i}{i!}.$$The {\em derangement graph} $\Gamma_n$
is the Cayley graph $\Gamma (S_n, D_n)$. That
is, two vertices $g, h$ of $\Gamma_n$ are joined if and only if $g(i ) \not= h(i )$ for all $i \in [n]$, or
equivalently $h^{-1}g$ fixes no point. Note that $\Gamma_n$ is loop-less because $D_n$ does not contain the identity element of $S_n$ and is a simple graph because $D_n$ is inverse-closed, that is, $D_n=\{g^{-1}~|~g\in D_n\}$. Clearly, $\Gamma_n$ is vertex-transitive, so it is $|D_n|$-regular. $\Gamma_n$ is connected if $n\ge 4$ because every vertex of $\Gamma_n$ can
be reached from the identity.

In the past few decades, there are lots of literatures investigate the edge-pancyclicicty and edge-fault-tolerant pancyclicicty of the Cayley graphs on symmetric group. Jwo et al.\cite{Jwo} and Tseng et al.\cite{Tseng} studied the bipancyclicity and edge-fault-tolerant bipancyclicicty of star graphs, respectively. Kikuchi and Araki \cite{Kikuchi} discussed the edge-bipancyclicity and edge-fault-tolerant bipancyclicicty of bubble-sort graphs. Tanaka et al.\cite{Tanaka} studied the bipancyclicity of Cayley graphs generated by transpositions. The derangement graph is also a well-studied Cayley graph on symmetric group. Over the years, there are many results about  the derangement graphs  such as the independence number \cite{KLW, M1}, EKR property \cite{MS}, eigenvalues \cite{KW18, KW13, KW10, R07}, automorphism group \cite{DZ11} and some other properties. Another important research content of the derangement graphs is the Hamiltonian property. The question as to whether the derangements graph is Hamiltonian was posed in \cite{R84, W88}.  Existence of a Hamilton cycle was shown in \cite{M85, W89}. In \cite{RS}, Rasmussen and Savage showed that the derangements graph is  Hamilton-connected,
that is, every pair of distinct vertices is joined by a Hamilton path. In this note, we will show that the derangements graph is edge pancyclic. The main result of this note is the following theorem.
\vskip.2cm

{\bf Theorem 1} {\em
 The derangements graph $\Gamma_n$ is edge pancyclic for $n\ge 4$.}

Then we have the following corollary directly.

 {\bf Corollary 2} {\em
 The derangements graph $\Gamma_n$ is (vertex) pancyclic for $n\ge 4$.}

\vskip.2cm

\section{Proof of Theorem 1}
\vskip.2cm

In order to proof our main result, we need the following lemma.

\vskip.2cm
{\bf Lemma 3 }[10, Theorem 45] {\em Let $\Gamma$ be a graph of order $n \ge 3$. If
$\delta(\Gamma)\ge (n+2)/2$, then $\Gamma$ is edge pancyclic.}
\vskip.2cm

We  need some extra notations. Let $S_n$ be the symmetric group on $[n] = \{1,\ldots, n\}$. We denote by $C_n$ the set of permutations in $S_n$ that consist of one single cycle of
length $n$. We call these {\em cyclic permutations}. It is clear that $|C_n| = (n -1)!$ and $\{1,\sigma(1),\sigma^2(1),\ldots,\sigma^{n-1}(1)\}=[n]$ for $\sigma\in C_n$. For $\sigma_1,\sigma_2\in S_n$, let $\Delta(\sigma_1,\sigma_2)$ be the numbers of the fixed points of $\sigma_1^{-1}\sigma_2$. We first have the following claim.

{\bf Claim 1} {\em For any $\alpha,\beta\in S_n$ and $\sigma\in C_n$, we have $$\sum_{i=0}^{n-1}\Delta(\alpha,\sigma^i\beta)=n.$$}

{\bf Proof of Claim 1} Note that for any $a,b\in [n]$, there is only $i\in \{0,1,\ldots,n-1\}$ such that $\sigma^i\beta(a)=b$. Since $\sigma\in C_n$, $\{\sigma^i\beta(a)~|~i=0,1,\ldots,n-1\}=[n]$ which implies the result holds.\q

\vskip.2cm

{\bf Proof of Theorem 1} Given $\alpha\beta\in E(\Gamma_n)$. By the definition of $\Gamma_n$, we have $\alpha\not=\beta$ and $ \Delta(\alpha,\beta)=0$. Since $|C_n|=(n-1)!$ and $n\ge 4$, there is $\sigma\in C_n$ such that $\alpha\beta^{-1}\not\in\{\sigma,\sigma^2,\ldots,\sigma^{n-1}\}$. By Claim 1 and $ \Delta(\alpha,\beta)=0$, there is $i_0\in [n-1]$ such that $ \Delta(\alpha,\sigma^{i_0}\beta)\ge 2$. Since $\alpha\beta^{-1}\not\in\{\sigma,\sigma^2,\ldots,\sigma^{n-1}\}$, $\alpha\not=\sigma^{i_0}\beta$. Let $\beta_0=\sigma^{i_0}\beta$ for short. Then there are $a,b,c,d\in [n]$ such that $\alpha(a)=\beta_0(a)=b$ and $\alpha(c)=\beta_0(c)=d$. Assume, without loss of generality, that $c=n$.

{\bf Claim 2} {\em We can assume that $d=n$.}

{\bf Proof of Claim 2} Assume $d\not=n$.  Let $\gamma$ be a transposition $(d,n)$ in $S_n$.
Denote a mapping $\varphi~:~V(\Gamma_n)\rightarrow V(\Gamma_n)$ such that $\varphi(\theta)=\gamma\theta$. Since $\alpha^{-1}\beta=\alpha^{-1}\gamma^{-1}\gamma\beta=(\gamma\alpha)^{-1}(\gamma\beta)$, $\varphi$ is an automorphism of $\Gamma_n$ which implies Claim 2 holds.\QEDopen

By Claim 2, we assume $\alpha(n)=\beta_0(n)=n$. Denote $T=\{\tau\in S_n~|~\tau(n)=n\}$. Then $|T|=(n-1)!$ and $\alpha,\beta_0\in T$. For any $\tau\in T$, let $A_\tau=\{\tau,\sigma\tau,\sigma^2\tau,\ldots,\sigma^{n-1}\tau\}$. Then $\beta\in A_{\beta_0}$.

{\bf Claim 3} {\em $S_n=\cup_{\tau\in T}A_\tau$ and $\Gamma_n[A_\tau]\cong K_n$ for any $\tau\in T$, where $K_n$ is a complete graph of order $n$.}

{\bf Proof of Claim 3} In order to show $S_n=\cup_{\tau\in T}A_\tau$, we just need to prove $A_{\tau_1}\cap A_{\tau_2}=\0$ for any $\tau_1,\tau_2\in T$ with $\tau_1\not=\tau_2$. Suppose there are $\tau_1,\tau_2\in T$ with $\tau_1\not=\tau_2$ such that $A_{\tau_1}\cap A_{\tau_2}\not=\0$. Then there are $i,j\in \{0,1,\ldots,n-1\}$ such that $\sigma^i\tau_1=\sigma^j\tau_2$. Assume $i>j$. Then we have $\sigma^{i-j}\tau_1=\tau_2$. Since $\tau_1,\tau_2\in T$, we have $\tau_1(n)=\tau_2(n)=n$ which implies $\sigma^{i-j}(n)=n$, a contradiction with $\sigma\in C_n$.

Let $\tau\in T$ and $\pi_1,\pi_2\in A_\tau$. Then there are $i,j\in \{0,1,\ldots,n-1\}$ such that $\pi_1=\sigma^i\tau$ and $\pi_2=\sigma^j\tau$. Assume $i>j$. If $\Delta(\pi_1,\pi_2)\ge 1$, say $\pi_1(k)=\pi_2(k)$ ($k\in [n]$), then $\sigma^{i-j}(\tau(k))=\tau(k)$, a contradiction with $\sigma\in C_n$. Hence $\Delta(\pi_1,\pi_2)=0$ and then $\Gamma_n[A_\tau]\cong K_n$.\QEDopen

{\bf Claim 4} {\em Let $\overline{\Gamma_{n-1}}$ be the complement of $\Gamma_{n-1}$. If $n\ge 5$, then $\overline{\Gamma_{n-1}}$ is edge pancyclic. If $n=4$, then $\overline{\Gamma_{3}}$ is edge even-pancyclic. }

{\bf Proof of Claim 4} Note that $\Gamma_{n-1}$ is $|D_{n-1}|$-regular. Then
$$\delta(\overline{\Gamma_{n-1}})=(n-1)!-1-(n-1)!\sum_{i=0}^{n-1}\frac{(-1)^i}{i!}=(n-1)!\sum_{i=1}^{n-1}\frac{(-1)^{i-1}}{i!}-1\ge \frac{(n-1)!}{2}+1$$ if $n\ge 5$. Thus $\overline{\Gamma_{n-1}}$ is edge pancyclic by Lemma 3 when $n\ge 5$.

If $n=4$, then $\overline{\Gamma_{3}}\cong K_{3,3}$. Thus the result holds.
\QEDopen

We complete the proof by considering the following two cases.

{\bf Case 1} $n\ge 5$.

Let $\tau=\tau(1)\tau(2)\cdots\tau(n-1)\tau(n)\in T$. Then $\tau(n)=n$. Denote $\widehat{\tau}=\tau(1)\cdots\tau(n-1)$ and $\widehat{T}=\{\widehat{\tau}~|~\tau\in T\}$. Then $\widehat{\tau}\in S_{n-1}$ and $\widehat{T}=S_{n-1}$. So $\Gamma(\widehat{T},D_{n-1})=\Gamma_{n-1}$.
Since $\alpha,\beta_0\in T$ and $ \Delta(\alpha,\beta_0)\ge 2$, $ \Delta(\widehat{\alpha},\widehat{\beta_0})\ge 1$ which implies $\widehat{\alpha}\widehat{\beta_0}\in E(\overline{\Gamma(\widehat{T},D_{n-1})})$. By Claim 4, for any integer $3\le k\le (n-1)!$, there are $\widehat{\tau_1},\widehat{\tau_2},\ldots,\widehat{\tau_k}\in \widehat{T}$ such that $\widehat{\tau_1}=\widehat{\alpha}$, $\widehat{\tau_2}=\widehat{\beta_0}$ and $\widehat{\tau_1}\widehat{\tau_2}\ldots\widehat{\tau_k}\widehat{\tau_1}$ is a cycle of $\overline{\Gamma(\widehat{T},D_{n-1})}$.
Since $\widehat{\tau_i}\widehat{\tau_{i+1}}\in E(\overline{\Gamma(\widehat{T},D_{n-1})})$,
$\Delta(\widehat{\tau_i},\widehat{\tau_{i+1}})\ge 1$ for all $1\le i\le k$, where the subscripts are modulo $k$. Hence $\Delta(\tau_i,\tau_{i+1})\ge 2$ for all $1\le i\le k$. By Claim 1, there is $\theta_{i+1}\in A_{\tau_{i+1}}$ such that $\Delta(\tau_i,\theta_{i+1})=0$ which implies $\tau_i\theta_{i+1}\in E(\Gamma_n)$ for all $1\le i\le k$. Recall $\tau_1=\alpha$, $\tau_2=\beta_0$, $\beta\in A_{\beta_0}$ and $ \Delta(\alpha,\beta)=0$. Then we can let $\theta_2=\beta$. By Claim 3, $\Gamma_n[A_{\tau_i}]\cong K_n$ for all $1\le i\le k$ and they are vertex-disjoint. Let $P_{ij}$ be a path of length $j$ connecting $\theta_i$ and $\tau_i$ in $\Gamma_n[A_{\tau_i}]$, where $1\le i\le k$ and $1\le j\le n-1$. Then
$$\theta_1P_{1j_1}\tau_1(=\alpha)\theta_2(=\beta)P_{2j_2}\tau_2\theta_3P_{3j_3}\ldots \theta_kP_{kj_k}\tau_k\theta_1$$is a cycle of length $k+\sum_{s=1}^kj_s$ contained the edge $\alpha\beta$, where $1\le j_s\le n-1$ for all $1\le s\le k$. Since $3\le k\le (n-1)!$, there is a cycle of length $l$ contained  $\alpha\beta$ for all $6\le l\le n!$. To finish our proof, we just need to show that there is a cycle of length $l$ contained  $\alpha\beta$ for all $3\le l\le 5$.

For any $\pi\in S_n$, denote $M(\pi)=\{(1,\pi(1)),(2,\pi(2)),\ldots,(n,\pi(n))\}$. Then there is a bijection between $\pi$ and $M(\pi)$. For any $\pi_1,\pi_2\in V(\Gamma_n)$, if $\pi_1\pi_2\in E(\Gamma_n)$, then $M(\pi_1)\cap M(\pi_2)=\0$; vice versa. Particularly, $M(\alpha)\cap M(\beta)=\0$.  We consider the complete bipartite graph $K_{n,n}$. Then $M(\pi)$ can be treated as a perfect matching of $K_{n,n}$. Since $n\ge 5$, we can find five disjoint perfect matchings $M(\alpha), M(\beta),M(\pi_1),M(\pi_2),M(\pi_3)$. Hence $\alpha\beta\pi_1\alpha$, $\alpha\beta\pi_1\pi_2\alpha$ and $\alpha\beta\pi_1\pi_2\pi_3\alpha$ are three cycles contained  $\alpha\beta$ of length 3, 4 and 5, respectively.

\begin{figure}[!htbp]
\begin{center}\includegraphics[scale=0.6]{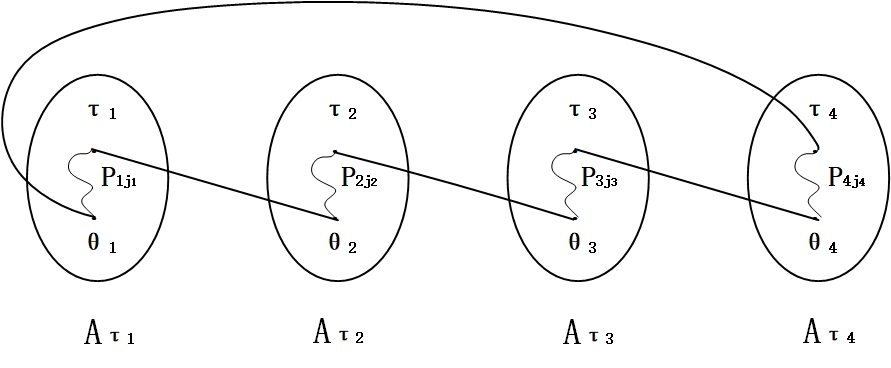}\\\end{center}
\caption{The construction of the cycle when $k=4$}\label{k4}
\end{figure}

{\bf Case 2} $n=4$.

Let $\tau=\tau(1)\tau(2)\tau(3)\tau(4)\in T$. Then $\tau(4)=4$. Denote $\widehat{\tau}=\tau(1)\tau(2)\tau(3)$ and $\widehat{T}=\{\widehat{\tau}~|~\tau\in T\}$. By the same argument as that of Case 1, we have $\Gamma(\widehat{T},D_{3})=\Gamma_{3}$. By Claim 4, for  $ k=4$ and $6$, there are $\widehat{\tau_1},\widehat{\tau_2},\ldots,\widehat{\tau_k}\in \widehat{T}$ such that $\widehat{\tau_1}=\widehat{\alpha}$, $\widehat{\tau_2}=\widehat{\beta_0}$ and $\widehat{\tau_1}\widehat{\tau_2}\ldots\widehat{\tau_k}\widehat{\tau_1}$ is a cycle of length $k$ in $\overline{\Gamma(\widehat{T},D_{3})}$. Then there is a cycle of length $l$ contained  $\alpha\beta$ for all $8\le l\le 24$. By the same argument, we can find four disjoint perfect matchings $M(\alpha), M(\beta),M(\pi_1),M(\pi_2)$. Hence $\alpha\beta\pi_1\alpha$ and $\alpha\beta\pi_1\pi_2\alpha$ are two cycles contained  $\alpha\beta$ of length 3 and 4, respectively. Now we consider $A_\alpha$ and $A_{\beta_0}$. Then $\beta\in A_{\beta_0}$. Recall $\beta_0=\sigma^{i_0}\beta$. Since $\Delta(\alpha,\beta)=0$, we have $\Delta(\sigma^{i_0}\alpha,\sigma^{i_0}\beta)=0$ which implies $\alpha_0\beta_0\in \Gamma_4$, where $\alpha_0=\sigma^{i_0}\alpha$. Since $\alpha_0\in A_\alpha$ and $|A_\alpha|=|A_{\beta_0}|=4$, we easily have cycles of length 5 to 7 contained  $\alpha\beta$ by Claim 3.

Thus we complete the proof.\q

\vskip.2cm

\end{document}